\documentclass[preprint,12pt]{article}

\usepackage{amssymb}

\newtheorem{thm}{Theorem}
\newtheorem{lemma}{Lemma}
\newtheorem{proposition}{Proposition}
\newtheorem{corollary}{Corollary}
\newtheorem{problem}{Problem}

\title{On the Number of Disjoint Pairs of S-permutation Matrices}
\author{Krasimir Yordzhev}
\date{}

\begin{document}
\unitlength=0.7mm \linethickness{0.5pt}
\maketitle

\begin{center}
Faculty of Mathematics and Natural Sciences\\
South-West University, Blagoevgrad, Bulgaria\\
e-mail: yordzhev@swu.bg
\end{center}

\begin{abstract}
 In [Journal of Statistical Planning and Inference (141) (2011) 3697--3704], Ro\-ber\-to Fontana offers an algorithm for obtaining Sudoku matrices. Introduced by Geir Dahl  concept disjoint pairs of S-permutation matrices [Linear Algebra and its Applications (430)  (2009) 2457–-2463] is used in this algorithm. Analyzing the works of G. Dahl and R. Fontana, the question of finding a general formula for counting disjoint pairs of $n^2 \times n^2$ S-permutation matrices  as a function of the integer $n$ naturally arises. This is an interesting combinatorial problem that deserves its consideration. The present work solves this problem. To do that, the graph theory techniques have been used. It has been shown that to count the number of disjoint pairs of $n^2 \times n^2$ S-permutation matrices, it is sufficient to obtain some numerical characteristics of the set of all  bipartite graphs of the type $g=\langle R_g \cup C_g , E_g \rangle$, where  $V=R_g \cup C_g$ is the set of vertices, and $E_g$ is the set of edges of the graph $g$, $R_g \cap C_g =\emptyset$, $|R_g |=|C_g |=n$.
\end{abstract}

\textbf{Keywords}: Binary matrix; S-permutation matrix; Sudoku matrix; Disjoint matrices; Bipartite graph

\textbf{2010 MSC}: 15B34; 05B20; 05C50

\section{Introduction}

Let $P_{ij}$, $1\le i,j\le n$ be $n^2$  square $n\times n$ matrices, which elements belong to the set $\mathfrak{Y}_{n^2} =\{ 1,2,\ldots ,n^2 \}$. Then $n^2 \times n^2$ matrix

\begin{equation}\label{P}
P =
\left[
\begin{array}{cccc}
P_{11} & P_{12} & \cdots & P_{1n} \\
P_{21} & P_{22} & \cdots & P_{2n} \\
\vdots & \vdots & \ddots & \vdots \\
P_{n1} & P_{n2} & \cdots & P_{nn}
\end{array}
\right]
\end{equation}
is called \emph{Sudoku Matrix}, if every row, every column and every submatrix $P_{ij}$, $1\le i,j\le n$ make permutation of the elements of set $\mathfrak{Y}_{n^2}$, i.e. every number $s\in \{ 1,2,\ldots ,n^2 \}$ is found just once in every row, every column and every submatrix $P_{ij}$. Submatrices $P_{ij}$ are called blocks of $P$.

Sudoku is a very popular game. On the other hand, it is well known
that Sudoku matrices are special cases of Latin squares in the class of gerechte designs \cite{Bailey}.

A matrix is called \emph{binary} if all of its elements are equal to 0 or 1. A square binary matrix is called \emph{permutation matrix}, if in every row and every column there is just one 1.

Let denote by $\Sigma_{n^2}$ the set of all permutation $n^2 \times n^2$ matrices of the following type
\begin{equation}\label{matrA}
A =
\left[
\begin{array}{cccc}
A_{11} & A_{12} & \cdots & A_{1n} \\
A_{21} & A_{22} & \cdots & A_{2n} \\
\vdots & \vdots & \ddots & \vdots \\
A_{n1} & A_{n2} & \cdots & A_{nn}
\end{array}
\right] ,
\end{equation}
where for every $s,t\in \{ 1,2,\ldots ,n\}$  $A_{st}$ is a  square $n\times n$ binary submatrix (block) with only one element equal to 1.
The elements of $\Sigma_{n^2} $ will be called \emph{S-permutation matrix}.

Two S-permutation  matrices $A$ and $B$ will be called \emph{disjoint}, if there are not $i,j\in \mathfrak{Y}_{n^2} =\{ 1,2,\ldots ,n^2 \}$ such that for the elements $a_{ij} \in A$ and $b_{ij}\in B$ the condition $a_{ij} =b_{ij} =1$ is satisfied.

The following obvious proposition is given in \cite{dahl}:

\begin{proposition}\label{disj} \cite{dahl}
A square $n^2 \times n^2$ matrix $P$  is Sudoku matrix if and only if there are $n^2$ mutually disjoint S-permutation  matrices $A_1 ,A_2 ,\ldots ,A_{n^2}$,  such that $P$ can be given in the following way:
$$P=1\cdot A_1 +2\cdot A_2 +\cdots +n^2 \cdot A_{n^2}$$
\end{proposition}
\hfill $\Box$

In \cite{Fontana} Roberto Fontana offers an algorithm which randomly gets a family of $n^2 \times n^2$ mutually disjoint S-permutation matrices, where $n=2,3$. In $n=3$ he run the algorithm 1000 times and found 105 different families of nine mutually disjoint S-permutation matrices. Then according to proposition \ref{disj} he decided that there are at least $9!\cdot 105  =38\; 102\; 400$ Sudoku matrices. This number is very small compared with the exact number of $9\times 9$ Sudoku matrices. In  \cite{Felgenhauer} it has been shown that  there are exactly
$$9! \cdot 72^2 \cdot 2^7 \cdot 27\; 704\; 267\; 971 = 6\; 670\; 903\; 752\; 021\; 072\; 936\; 960 $$
number of $9\times 9$ Sudoku matrices.

To evaluate the effectiveness of Fontana's algorithm is necessary to calculate the probability of two randomly generated matrices to be disjoint.
As it is proved in \cite{dahl} the number of S-permutation matrices is equal to $\displaystyle \left( n! \right)^{2n} $.
Thus the question of enumerating all disjoint pairs of  S-permutation matrices naturally arises. This work is devoted to this task.
In  the  present  paper, we will formulate and prove  a  formula  (Theorem  \ref{t1} and in another form Theorem \ref{t2}) for counting  disjoint pairs  of   S-permutation matrices. To do that, we will use graph theory techniques.

\section{A  formula  for counting  disjoint pairs  of   S-permutation matrices--the graph theory approach}

Let us denote by $\mathcal{S}_m$ the symmetric group of order $m$, i.e. the group  of all one-to-one mapping of the  set $\mathfrak{Y}_{m} =\left\{ 1,2,\ldots ,m\right\}$ in itself. If $x\in \mathfrak{Y}_m $, $\rho\in \mathcal{S}_m$, then the image of the element $x$ in the mapping $\rho$ we denote by $\rho (x)$.

\emph{Bipartite graph} is the ordered triplet
 $$
g=\langle R_g \cup C_g , E_g \rangle ,
$$
where $R_g$ and $C_g$ are non-empty sets such that $R_g \cap C_g =\emptyset$, the elements of which will be called \emph{vertices}. $E_g \subseteq R_g \times C_g =\{ \langle r,c \rangle \; |\; r\in R_g ,c\in C_g \}$  - the set of \emph{edges}. Repeated edges are not allowed in our considerations.

Let $g'=\langle R_{g'} \cup C_{g'} , E_{g'} \rangle$ and $g''=\langle R_{g''} \cup C_{g''} , E_{g''} \rangle$. We will say that the graphs $g'$ and $g''$ are \emph{isomorphic} and we will write $g'\cong g''$, if $R_{g'} =R_{g''}$, $C_{g'} =C_{g''}$ and there are $\rho \in \mathcal{S}_m$ and $\sigma \in\mathcal{S}_n$, where $m=|R_{g'}|=|R_{g''}|$ and $n=|C_{g'}|=|C_{g''}|$, such that $\langle r, c\rangle \in E_{g'} \Longleftrightarrow \langle \rho (r) ,\sigma (c)\rangle \in E_{g''}$. The object of this work is bipartite graphs considered to within isomorphism.

For more details on graph theory see  for example \cite{diestel} or  \cite{harary}.

Let $n$ and $k$ be two positive integers and let $0\le k\le n^2$. Let us denote by $\mathfrak{G}_{n,k} $ the set of all  bipartite graphs without multiple edges of the type $g=\langle  R_g \cup C_g ,E_g \rangle $, considered to within isomorphism, such that $ |R_g |= |C_g |= n$ and $|E_g |=k$.

Let $g=\langle R_g \cup C_g ,E_g \rangle \in \mathfrak{G}_{n,k}$ for some natural numbers $n$ and $k$ and let $v\in V_g =R_g \cup C_g$. By $\gamma (v)$ we denote the set of all vertices of $V_g$, incident with $v$, i.e. $u\in \gamma (v)$ if and only if there is an edge in $E_g$ connecting $u$ and $v$. If $v$ is an isolated vertex (i.e. there is no edge, incident with $v$), then by definition $\gamma (v)=\emptyset$ and $|\gamma (v)|=0$. If $v\in R_g$, then obviously $\gamma (v)\subseteq C_g$, and if $v\in C_g$, then $\gamma (v)\subseteq R_g$. Obviously
$$\sum_{v\in V_{g}} |\gamma (v)|=2k.$$

Let $g=\langle R_g \cup C_g ,E_g \rangle \in \mathfrak{G}_{n,k}$ and let $u,v\in V_g =R_g\cup C_g$. We will say that $u$ and $v$ are equivalent and we will write $u\sim v$ if $\gamma (u) =\gamma (v)$. If $u$ and $v$ are isolated, then by definition $u\sim v$ if and only if $u$ and $v$ belong simultaneously to $R_g$, or $C_g$. The above introduced relation is obviously an equivalence relation.

By ${V_g}_{/\sim} $ we denote the obtained factor-set (the set of the equivalence classes) according to relation $\sim$ and let
$${V_g}_{/\sim} =\left\{ \Delta_1 ,\Delta_2 ,\ldots ,\Delta_s \right\} ,$$
where $\Delta_i \subseteq R_g$, or $\Delta_i \subseteq C_g$, $i=1,2,\ldots s$, $2\le s \le 2n$. We put
$$\delta_i =|\Delta_i |,\quad 1\le \delta_i \le n , \quad i=1,2,\ldots , s$$
and for every $g\in \mathfrak{G}_{n,k}$ we define multi-set (set with repetition)
$$\left[ g \right] =\left\{ \delta_1 ,\delta_2 ,\ldots \delta_s \right\} ,$$
where $\delta_1 ,\delta_2 ,\ldots ,\delta_s$ are natural numbers, obtained by the above described way.
 Obviously
 $$\sum_{i=1}^s \delta_i =2n .$$

If $z_1 \; z_2 \; \ldots \; z_n$ is a permutation of the elements of the set $\mathfrak{Y}_n =\left\{ 1,2,\ldots ,n\right\}$ and we shortly denote $\rho$ this permutation, then in this case we denote by $\rho (i )$ the $i$-th element of this permutation, i.e. $\rho (i) =z_i$, $i=1,2,\ldots ,n$.

\begin{lemma}\label{l1}
The number of all $n\times n$ binary matrices with exactly $k$, $k=0,1,\ldots ,n^2$ ones, denoted by $b(n,k)$ is equal to
\begin{equation}\label{bnk}
b(n,k)=(n!)^2 \sum_{g\in \mathfrak{G}_{n,k} } \frac{1}{\displaystyle \prod_{\delta \in [g]} \delta !}
\end{equation}
\end{lemma}

Proof. Let $A=\left[ a_{ij} \right]_{n\times n}$ be $n\times n$ binary matrix with exactly $k$ ones. Then we construct graph $g=\langle R_g\cup C_g ,E_g\rangle $, such that the set $R_g=\{ r_1 ,r_2 ,\ldots ,r_n \}$  corresponds to the rows of $A$, and $C_g=\{ c_1 , c_2 ,\ldots , c_n \}$  corresponds to the columns of $A$, however there is an edge connecting the vertices $r_i$ and $c_j$ if and only if $a_{ij} =1$. The graph, which has been constructed, obviously belongs to $\mathfrak{G}_{n,k}$.

Conversely, let $g=\langle R_g\cup C_g ,E_g\rangle \in \mathfrak{G}_{n,k}$. We number at a random way the vertices of $R_g$ with the help of the natural numbers from 1 to $n$ without repeating any of the numbers. We analogously number the vertices of $C_g$. This can be made by $(n!)^2$ ways. Then we construct binary $n\times n$ matrix $A=\left[a_{ij} \right]_{n\times n}$, such that $a_{ij} =1$ if and only if there is an edge in $E_g$ connecting the vertex with number $i$ of $R_g$ with the vertex with number $j$ of $C_g$. As $g\in \mathfrak{G}_{n,k}$, then the matrix, that has been constructed, has exactly $k$ units. It is easy to see that when $q,r\in \mathfrak{Y}_n$, $q$-th  and $r$-th rows of $A$ are equal to each other (i.e. the matrix $A$ does not change if we changes the places of these two rows) if and only if the vertices of $R_g$ corresponding to numbers $q$ and $r$ are equivalent according to relation $\sim$. Analogous assertion is true about the columns of the matrix $A$ and the edges of the set $C_g$, which proves formula (\ref{bnk}).

\hfill $\Box$

\begin{corollary}\label{lrl1}
For every positive integer $n$, $n\ge 2$ and every natural number $k$ such that $0\le k\le n^2$ the following formula is true:
\begin{equation}\label{form5}
(n!)^2 \sum_{g\in \mathfrak{G}_{n,k} } \frac{1}{\displaystyle \prod_{\delta \in [g]} \delta !} ={n^2 \choose k}
\end{equation}
\hfill $\Box$
\end{corollary}

Let $\Pi_n$ denotes the set of all $n\times n$ matrices, constructed such that $\pi\in\Pi_n$ if and only if the following three conditions are true:

i) the elements of $\pi$ are ordered pairs of numbers $\langle i,j\rangle$, where $1\le i,j\le n$;

ii) if
$$\left[ \langle a_1 \quad b_1 \rangle \quad \langle a_2 ,b_2 \rangle \quad \cdots  \quad \langle a_n ,b_n \rangle \right]$$
is the $i$-th row of $\pi$ for any $i\in \mathfrak{Y}_n =\{ 1,2,\ldots ,n\}$, then the numbers $a_1 \; a_2 \; \ldots \; a_n$ in this order construct permutation of the elements of the set $\mathfrak{Y}_n$

iii) if
$$\left[
\begin{array}{c}
\langle a_1 ,b_1 \rangle \\
\langle a_2 ,b_2 \rangle \\
\vdots \\
\langle a_n ,b_n \rangle \\
\end{array}
\right]
$$
is the $j$-th column of $\pi$ for any $j\in \mathfrak{Y}_n$, then numbers $b_1 ,b_2 ,\ldots , b_n$ in this order construct permutation of the elements of the set $\mathfrak{Y}_n$.

From the definition it follows that every row and every column of any matrix of the set $\Pi_n$ can be identified with permutation of elements of the set $\mathfrak{Y}_n$. Conversely for every $(2n)$-tuple $\langle \langle \rho_1 ,\rho_2 ,\ldots ,\rho_n \rangle ,\langle \sigma_1 ,\sigma_2 ,\ldots , \sigma_n \rangle \rangle$, where $\rho_i = \rho_i (1)\; \rho_i (2) \; \ldots \; \rho_i (n)$, $\sigma_j = \sigma_j (1)\; \sigma_j (2)\; \ldots \; \sigma_j (n)$, $1\le i,j\le n$ are permutations of elements of $\mathfrak{Y}_n$, then the matrix
$$
\pi =
\left[
\begin{array}{cccc}
\langle \rho_1 (1),\sigma_1 (1)\rangle & \langle \rho_1 (2),\sigma_2 (1)\rangle & \cdots & \langle \rho_1 (n),\sigma_n (1)\rangle \\
\langle \rho_2 (1),\sigma_1 (2)\rangle & \langle \rho_2 (2),\sigma_2 (2)\rangle & \cdots & \langle \rho_2 (n),\sigma_n (2)\rangle \\
\vdots & \vdots & \ddots & \vdots \\
\langle \rho_n (1),\sigma_1 (n)\rangle  & \langle \rho_n (2),\sigma_2 (n)\rangle & \cdots & \langle \rho_n (n),\sigma_n (n)\rangle
\end{array}
\right]
$$
is matrix of $\Pi_n$. Hence
\begin{equation}\label{|Pin|}
\left| \Pi_n \right| =\left( n! \right)^{2n}
\end{equation}

We say that matrices $\pi ' ,\pi '' \in\Pi_n$, where $\pi ' =\left[ {p'}_{ij} \right]_{n\times n}$, $\pi '' =\left[ {p''}_{ij} \right]_{n\times n}$ are \textit{disjoint}, if ${p'}_{ij} \ne {p''}_{ij}$ for every pair of indices $i,j\in\mathfrak{Y}_n$.

\begin{lemma}\label{l2}
Let $n$ be a positive integer, $n\ge 2$. Then there is one to one correspondence between the sets $\Sigma_{n^2}$ and $\Pi_n$ .
\end{lemma}

Proof. Let $A\in \Sigma_{n^2}$. Then $A$ is constructed with the help of formula (\ref{matrA}) and for every $i,j\in \mathfrak{Y}_n$ in the block $A_{ij} $ there is only one 1 and let this 1 has coordinates $(a_i ,b_j )$. For every $i,j\in \mathfrak{Y}_n$ we obtain ordered pairs of numbers $\langle a_i ,b_j \rangle$ corresponding to these coordinates. As in every row and every column of $A$ there is only one 1, then the matrix $\left[ \alpha_{ij} \right]_{n\times n}$, where $\alpha_{ij} =\langle a_i ,b_j \rangle $, $1\le i,j\le n$, which is obtained by the ordered pairs of numbers is matrix of $\Pi_n$, i.e. matrix for which the conditions i), ii) and iii) are true.

Conversely, let $\left[ \alpha_{ij} \right]_{n\times n} \in \Pi_n$, where $\alpha_{ij} =\langle a_i ,b_j \rangle $, $i,j \in \mathfrak{Y}_n$, $a_i ,b_j \in \mathfrak{Y}_n$. Then for every $i,j\in \mathfrak{Y}_n$ we construct binary $n\times n$ matrices $A_{ij}$ with only one 1 with coordinates $(a_i ,b_j )$. Then we obtain the matrix of type: (\ref{matrA}). According to the properties i), ii) and iii), it is obvious that the obtained matrix is S-permutation matrix.

\hfill $\Box$

From Lemma \ref{l2} and formula (\ref{|Pin|}) it follows the proof of Proposition \ref{disj} in \cite{dahl}.

\begin{corollary} {\rm (Proposition \ref{disj} of \cite{dahl})}
The number of S-permutation matrices of order $n^2$ is given by
\begin{equation}\label{fcrl2}
\left| \Sigma_{n^2} \right| = \left( n! \right)^{2n}
\end{equation}
\hfill $\Box$
\end{corollary}

It is easy to see that with respect of the described in Lemma \ref{l2} one to one correspondence, every pair of disjoint matrices of $\Sigma_{n^2}$ will correspond to a pair of disjoint matrices of $\Pi_n$ and conversely every pair of disjoint matrices of $\Pi_n$ will correspond to a pair of disjoint matrices of $\Sigma_{n^2}$.

\begin{corollary}
The number of all  pairs of disjoint matrices of $\Sigma_{n^2}$ is equal to the number of all pairs of disjoint matrices of $\Pi_n$.

\hfill $\Box$
\end{corollary}

\begin{lemma}\label{l3}
The number $q(n,k)$ of all  ordered pairs of matrices $\langle \pi' ,\pi'' \rangle$, where $\pi' ,\pi'' \in\Pi_n$ and  having  at least $k$, $k=0,1,\ldots ,n^2$ equal elements is equal to
\begin{equation}\label{fl3}
q(n,k)=(n!)^{2(n+1)} \sum_{g\in \mathfrak{G}_{n,k} } \frac{\displaystyle \prod_{v\in R_g \cup C_g} ( n-|\gamma (v)|)!}{\displaystyle \prod_{\delta \in [g]} \delta !}
\end{equation}

\end{lemma}

Proof. Let $\pi' =\left[ p'_{ij} \right]_{n\times n} ,\pi'' =\left[ p''_{ij} \right]_{n\times n} \in \Pi_n$ and let $\pi'$ and $\pi''$ have exactly $k$ equal elements (by components). Then we definitely obtain the binary $n\times n$ matrix $A=\left[ a_{ij} \right]_{n\times n}$, such that $a_{ij} =1$ if and only if $p'_{ij} =p''_{ij}$, $i,j\in\mathfrak{Y}_n$. According to Lemma \ref{l1} there is only one graph $g\in \mathfrak{G}_{n,k}$ corresponding to the matrix $A$ and respectively to the ordered pair of matrices $\langle \pi' ,\pi'' \rangle \in \Pi\times\Pi$.

Conversely, let $g=\langle R_g \cup C_g ,E_g \rangle \in \mathfrak{G}_{n,k}$ and let $V_g =R_g \cup C_g$. Then there are $\displaystyle \frac{\displaystyle (n!)^2}{\displaystyle \prod _{\delta\in [g]} \delta !}$ in number $n\times n$ binary matrices, corresponding to $g$, following the described in Lemma \ref{l1} rule and let $A=[a_{ij} ]_{n\times n}$ be one of them.
Let $\pi =\left[ p_{ij} \right]_{n\times n}$ be an arbitrary  matrix of $\Pi_n$. We will find the number $h(\pi ,A)$  of all matrices  $\pi' =[p'_{ij} ]_{n\times n} \in \Pi_n$, such that $p'_{ij} =p_{ij}$ if $a_{ij}=1$. It is allowed to exist $s,t\in \mathfrak{Y}_n$ such that $a_{st} =0$ and $p'_{st} =p_{st}$. Let the $i$-th row ($i=1,2,\ldots ,n$) of $\pi$ correspond to the permutation $\rho_i$ of the elements of $\mathfrak{Y}_n$ and let $i$-th row of the matrix $A$  correspond to the vertex  $r_i \in R_g$ of the graph $g$ (see Lemma \ref{l1}). Then there are $(n-|\gamma (r_i )|)!$ in number permutations $\rho'_i$ of elements of $\mathfrak{Y}_n$, such that if $a_{it} =1$, then $\rho_i (t) =\rho'_i (t)$, $t\in\mathfrak{Y}_n$. Analogously we can prove the corresponding assertion  for the columns  of $\pi$. Hence, $\displaystyle h(\pi ,A) =\prod_{v\in V_g} (n-|\gamma (v)|)!$. From all the things that have been described above, it follows that for every $\pi\in\Pi_n$ there are
$$\sum_{g\in\mathfrak{G}_{n,k}} \frac{\displaystyle (n!)^2}{\displaystyle \prod_{\delta\in [g]} \delta !} \prod_{v\in V_g} (n-|\gamma (v)|)!$$
matrices of $\Pi_n$, which have at least $k$ elements, which are equal to the corresponding elements of $\pi$ (by components). But according to formula (\ref{|Pin|}) $|\Pi_n |=(n!)^{2n}$, from where it follows formula (\ref{fl3}).

\hfill $\Box$

The purpose of the present work is to prove the following theorem:

\begin{thm}\label{t1}
Let $n\ge 2$ be a positive integer. Then the number $D_{n^2}$  of all possible ordered pairs of disjoint matrices in $\Sigma_{n^2}$ is equal to
\begin{equation}\label{f2}
D_{n^2} = (n!)^{4n} + (n!)^{2(n+1)}\sum_{k=1}^{n^2} (-1)^k   \sum_{g\in \mathfrak{G}_{n,k} } \frac{\displaystyle \prod_{v\in R_g \cup C_g} ( n-|\gamma (v)|)!}{\displaystyle \prod_{\delta \in [g]} \delta !}.
\end{equation}

The number $d_{n^2}$ of all non-ordered pairs of disjoint matrices in $\Sigma_{n^2}$ is equal to
\begin{equation}\label{f3}
d_{n^2} =\frac{1}{2} D_{n^2}
\end{equation}
\end{thm}

Proof. Let $n$ be a positive integer, $n\ge 2$. Then applying Lemma \ref{l2} and the inclusion-exclusion principle we obtain, that the number $D_{n^2}$ of all possible ordered pairs of disjoint matrices of $\Sigma_{n^2}$ is equal to
$$
D_{n^2} =   |\Pi_n |^2 +\sum_{k=1}^{n^2} (-1)^k q(n,k),
$$
where the function $q(n,k)$ is defined in Lemma \ref{l3} and may be calculated with the help of formula (\ref{fl3}), and $|\Pi_n |$ with the help of formula (\ref{|Pin|}). Then we obtain the proof of formula (\ref{f2}).

Having in mind that the ''disjoint'' relation is symmetric and irreflexive we obtain formula (\ref{f3}).

Theorem \ref{t1} is proved.

\hfill $\Box$

Let $n$ and $k$ are positive integers and let $g\in \mathfrak{G}_{n,k}$. We examine the ordered $(n+1)$-tuple

\begin{equation}\label{PSI}
\Psi (g)=\langle \psi_0 (g) ,\psi_1 (g),\ldots ,\psi_n (g)\rangle ,
\end{equation}
where $\psi_i (g)$, $i=0,1,\ldots ,n$ is equal to the number of vertices of $g$ incidents with exactly $i$ number of edges. It is obvious that $\displaystyle \sum_{i=1}^n i\psi_i (g)=2k$ is true for all $g\in \mathfrak{G}_{n,k}$. Then formula (\ref{f2}) can be presented
$$
D_{n^2} = (n!)^{4n} + (n!)^{2(n+1)}\sum_{k=1}^{n^2} (-1)^k   \sum_{g\in \mathfrak{G}_{n,k} } \frac{\displaystyle \prod_{i=0}^n \left[ \left( n-i\right) ! \right]^{\psi_i (g)}}{\displaystyle \prod_{\delta \in [g]} \delta !}.
$$

Since $(n-n)!=0!=1$ and $[n-(n-1)]!=1!=1 $, then

\begin{equation}\label{f2star}
D_{n^2} = (n!)^{4n} + (n!)^{2(n+1)}\sum_{k=1}^{n^2} (-1)^k   \sum_{g\in \mathfrak{G}_{n,k} } \frac{\displaystyle \prod_{i=0}^{n-2} \left[ \left( n-i\right) ! \right]^{\psi_i (g)}}{\displaystyle \prod_{\delta \in [g]} \delta !}.
\end{equation}

Consequently, to apply formula (\ref{f2star}) for each bipartite graph $g\in \mathfrak{G}_{n,k}$ and for the set $\mathfrak{G}_{n,k}$ of  bipartite graphs, it is necessary to obtain the following numerical characteristics:

\begin{equation}\label{star2}
\omega (g) = \frac{\displaystyle \prod_{i=0}^{n-2} \left[ \left( n-i\right) ! \right]^{\psi_i (g)}}{\displaystyle \prod_{\delta \in [g]} \delta !}
\end{equation}
and
\begin{equation}\label{thet}
\theta (n,k) =\sum_{g\in \mathfrak{G}_{n,k} } \omega (g)
\end{equation}

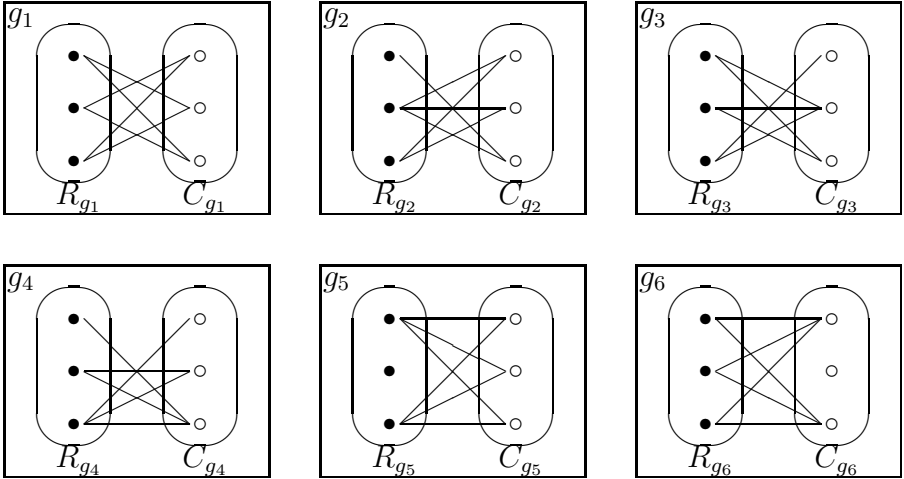
\begin{figure}[h]
\begin{picture}(170,90)

\multiput(0,0)(60,0){3}{
\multiput(0,0)(0,50){2}{

\put(0,0){\framebox(50,40)}

\multiput(13,21)(24,0){2}{\oval(14,30)}

\put(13,30){\circle*{2}}
\put(37,30){\circle{2}}

\put(13,20){\circle*{2}}
\put(37,20){\circle{2}}

\put(13,10){\circle*{2}}
\put(37,10){\circle{2}}
}
}

\put(3,87){\makebox(0,0){$g_{1}$}}
\put(63,87){\makebox(0,0){$g_{2}$}}
\put(123,87){\makebox(0,0){$g_{3}$}}

\put(14,53){\makebox(0,0){$R_{g_{1}}$}}
\put(38,53){\makebox(0,0){$C_{g_{1}}$}}

\put(74,53){\makebox(0,0){$R_{g_{2}}$}}
\put(98,53){\makebox(0,0){$C_{g_{2}}$}}

\put(134,53){\makebox(0,0){$R_{g_{3}}$}}
\put(158,53){\makebox(0,0){$C_{g_{3}}$}}

\put(3,37){\makebox(0,0){$g_{4}$}}
\put(63,37){\makebox(0,0){$g_{5}$}}
\put(123,37){\makebox(0,0){$g_{6}$}}

\put(14,3){\makebox(0,0){$R_{g_{4}}$}}
\put(38,3){\makebox(0,0){$C_{g_{4}}$}}

\put(74,3){\makebox(0,0){$R_{g_{5}}$}}
\put(98,3){\makebox(0,0){$C_{g_{5}}$}}

\put(134,3){\makebox(0,0){$R_{g_{6}}$}}
\put(158,3){\makebox(0,0){$C_{g_{6}}$}}

\put(15,80){\line(1,-1){20}}
\put(15,80){\line(2,-1){20}}

\put(15,70){\line(2,1){20}}
\put(15,70){\line(2,-1){20}}

\put(15,60){\line(1,1){20}}
\put(15,60){\line(2,1){20}}

\put(75,80){\line(1,-1){20}}

\put(75,70){\line(1,0){20}}
\put(75,70){\line(2,1){20}}
\put(75,70){\line(2,-1){20}}

\put(75,60){\line(1,1){20}}
\put(75,60){\line(2,1){20}}

\put(135,80){\line(1,-1){20}}
\put(135,80){\line(2,-1){20}}

\put(135,70){\line(1,0){20}}
\put(135,70){\line(2,-1){20}}

\put(135,60){\line(1,1){20}}
\put(135,60){\line(2,1){20}}

\put(15,30){\line(1,-1){20}}

\put(15,20){\line(1,0){20}}
\put(15,20){\line(2,-1){20}}

\put(15,10){\line(1,0){20}}
\put(15,10){\line(1,1){20}}
\put(15,10){\line(2,1){20}}

\put(75,30){\line(1,0){20}}
\put(75,30){\line(1,-1){20}}
\put(75,30){\line(2,-1){20}}

\put(75,10){\line(1,0){20}}
\put(75,10){\line(1,1){20}}
\put(75,10){\line(2,1){20}}

\put(135,30){\line(1,0){20}}
\put(135,30){\line(1,-1){20}}

\put(135,20){\line(2,1){20}}
\put(135,20){\line(2,-1){20}}

\put(135,10){\line(1,0){20}}
\put(135,10){\line(1,1){20}}

\end{picture}
\caption{The set of bipartite graphs $\mathfrak{G}_{3,6}$}\label{n3k6}
\end{figure}

\textbf{Example}. When $n = 3$ and $k = 6$ the set $\mathfrak{G}_{3,6}= \{ g_{1} ,g_{2} ,g_{3} ,g_{4} ,g_{5} ,g_{6} \}$ consists of six bipartite graphs, which are shown in Figure \ref{n3k6}.

For graph $g_{1} \in\mathfrak{G}_{3,6}$ we have:

$$\left[ g_{1} \right] =\left\{ 1,1,1,1,1,1\right\}$$

$$\Psi (g_{1} ) =\langle 0,0,6,0\rangle$$

$$\omega (g_{1} )=\frac{\left[ (3-0)!\right]^0 \left[ (3-1)!\right]^0 }{1!\; 1!\; 1!\; 1!\; 1!\; 1!}  =1$$

For graphs $g_{2} \in \mathfrak{G}_{3,6}$ and $g_{3} \in \mathfrak{G}_{3,6}$ we have:

$$\left[ g_{2} \right] =\left[ g_{3} \right] =\left\{ 1,1,1,1,2\right\}$$

$$\Psi (g_{2} ) =\Psi (g_{3} ) =\langle 0,1,4,1\rangle$$

$$\omega (g_{2} )=\omega (g_{3} )=\frac{\left[ (3-0)!\right]^0 \left[ (3-1)!\right]^1 }{1!\; 1!\; 1!\; 1!\; 2!} = \frac{6^0 \cdot 2^1}{2} =1$$

For graph $g_{4} \in\mathfrak{G}_{3,6}$ we have:

$$\left[ g_{4} \right] =\left\{ 1,1,1,1,1,1\right\}$$

$$\Psi (g_{4} ) =\langle 0,2,2,2\rangle$$

$$\omega (g_{4} )=\frac{\left[ (3-0)!\right]^0 \left[ (3-1)!\right]^2 }{1!\; 1!\; 1!\; 1!\; 1!\; 1!} =\frac{6^0 \cdot 2^2}{1}  =4$$

For graphs $g_{5} \in \mathfrak{G}_{3,6}$ and $g_{6} \in \mathfrak{G}_{3,6}$ we have:

$$\left[ g_{5} \right] =\left[ g_{6} \right] =\left\{ 1,2,3\right\}$$

$$\Psi (g_{5} ) =\Psi (g_{6} ) =\langle 1,0,3,2\rangle$$

$$\omega (g_{5} )=\omega (g_{6} )=\frac{\left[ (3-0)!\right]^1 \left[ (3-1)!\right]^0 }{1!\; 2!\; 3!} = \frac{6^1 \cdot 2^0}{2\cdot 6} =\frac{1}{2}$$

Then for the set $\mathfrak{G}_{3,6}$ we get:
$$
\theta (3,6)=\sum_{g\in \mathfrak{G}_{3,6}} \omega (g)=1+1+1+4+\frac{1}{2} +\frac{1}{2} =8
$$

\hfill $\Box$

All bipartite graphs with n = 2 and n = 3 together with the above numerical characteristics are described in \cite{Yordzhev}.

Using the numerical characteristics (\ref{star2}) and (\ref{thet}), we obtain the following variety  of Theorem \ref{t1}:

\begin{thm}\label{t2}
\begin{equation}\label{endstar}
D_{n^2} = (n!)^{4n} + (n!)^{2(n+1)}\sum_{k=1}^{n^2} (-1)^k   \theta (n,k),
\end{equation}
where $\theta (n,k)$ is described using formulas (\ref{thet}) and (\ref{star2}).

\hfill $\Box$
\end{thm}

\section{Conclusions and future work}

In order to apply Theorem \ref{t1} or Theorem \ref{t2} it is necessary to describe all bipartite graphs to within isomorphism $g=\langle R_g \cup C_g , E_g \rangle $, where $|R_g |=|C_g |=n$. For n = 2 and n = 3 the application of Theorem \ref{t2} is described in \cite{Yordzhev} in full details. According to calculations, we obtain:

\begin{itemize}
\item The number $D_4$ of all disjoint ordered pairs of  S-permutation  matrices in $ n =  2$ is equal to
\begin{equation}\label{D_4}
D_4 = 144
\end{equation}
\item The number $d_4$ of all disjoint non-ordered pairs of S-permutation  matrices in $ n =  2$ is equal to
\begin{equation}\label{ddd2}
d_4 =\frac{1}{2} D_4 =72
\end{equation}
\item The number $D_9$ of all disjoint ordered pairs of  S-permutation  matrices in $ n =  3$ is equal to
\begin{equation}\label{D_9}
D_9 = 1\; 260\; 085\; 248
\end{equation}
\item The number $d_4$ of all disjoint non-ordered pairs of S-permutation  matrices in $ n =  3$ is equal to
\begin{equation}\label{ddd3}
d_9 =\frac{1}{2} D_9 =630\; 042\; 624
\end{equation}
\end{itemize}

\begin{problem}\label{prbl1}
Let $n\ge 2$ is a natural number and let $G$ be a simple graph having $(n!)^{2n}$ vertices. Let each vertex of $G$ be identified with an element of the  set $\Sigma_{n^2}$ of all $n^2 \times n^2$ S-permutation matrices. Two vertices are connected by an edge if and only if the corresponding matrices are disjoint. The problem is to find the number of all complete subgraphs of $G$ having $n^2$ vertices:
\end{problem}

Note that the number of edges in graph $G$ is equal to $d_{n^2}$ and can be calculated using formula (\ref{f2}) and formula (\ref{f3})  (respectively formulas (\ref{star2}), (\ref{thet}), (\ref{endstar})  and (\ref{f3})).

Denote by $z_n$  the solution of the Problem \ref{prbl1} and let $\sigma_n$ is the number of all $n^2 \times n^2$ Sudoku matrices. Then according to Proposition \ref{disj} and the method of construction of the graph
 $G$, it follows that the next equality is valid:
\begin{equation}\label{z_n}
z_n =\frac{\sigma_n}{(n^2 )!}
\end{equation}

We do not know a general formula for finding the number of all $n^2 \times n^2$ Sudoku matrices for each natural number $n\ge 2$ and we consider that this is an open combinatorial problem. Only some special cases are known. For example in $n=2$ it is known that $\sigma_2 =288$ \cite{yorkost}. Then according to formula (\ref{z_n}) we get:

$$z_2 =\frac{\sigma_2}{4!} =\frac{288}{24} =12$$

In  \cite{Felgenhauer} it has been shown that in $n=3$ there are exactly
$$\sigma_3 = 6\; 670\; 903\; 752\; 021\; 072\; 936\; 960 =$$
$$= 9! \times 72^2 \times 2^7 \times 27\; 704\; 267\; 971 =$$
$$2^{20} \times  3^8 \times  5^1 \times  7^1 \times 27\; 704\; 267\; 971^1  \sim 6.671\times 10^{21}$$
number of Sudoku matrices. Then according to formula (\ref{z_n}) we get:

$$z_3 =\frac{\sigma_3}{9!} =\frac{6\; 670\; 903\; 752\; 021\; 072\; 936\; 960}{362\; 880} =18\; 383\; 222\; 420\; 692\; 992$$

\bibliographystyle{plain}
\bibliography{S-permMatr}

\end{document}